\documentclass[11pt]{article}
\newtheorem{thm}{Theorem}
\newtheorem{lem}{Lemma}
\begin{document}
\title{Another simple reformulation of the four color theorem}
\author{Ajit A. Diwan\\
Department of Computer Science and Engineering\\
Indian Institute of Technology Bombay\\
Powai, Mumbai 400076, India.}
\maketitle

\noindent {\bf Abstract:} We give a simple reformulation of the four color 
theorem as a problem on strings over a four letter alphabet. 

\section{Introduction}
The four color theorem is one of the cornerstones of graph theory. While there
are several equivalent statements and generalizations within graph theory, the
theorem has also been shown to be equivalent to statements involving mathematical
objects other than graphs. Some examples of these include arithmetic and
algebraic formulations~\cite{K,M} and also in terms of formal languages~\cite{CRZ,EL}.
We give a simple reformulation involving strings over a finite
alphabet. We show that the problem reduces to showing that there is no path
between two specific states in a specific finitely branching automaton with countably 
infinite states.  While we do not know any general techniques for doing this, we
hope that this formulation involving only strings may yield a simpler proof.

\section{Formulation}
Let $A =$ $\{$\texttt{a},\texttt{b},\texttt{c},\texttt{d}$\}$ be the
alphabet and let $A^*$ be the set of all finite length strings over $A$. A subset $L 
\subset A^*$ is called an $l$-subset if every string in $L$ has length
exactly $l$. Let $\mathcal{L}$ denote the collection of all $l$-subsets of $A^*$ for
all $l \ge 0$. Thus  a subset $X \subseteq A^*$ is in $\mathcal{L}$ iff
$X$ is an $l$-subset for some integer $l \ge 0$. We construct an automaton
whose set of states is $\mathcal{L}$.

Let $s = s_1s_2\ldots s_l$ be a string of length $l \ge 3$ over $A$. Let 
$1 \le i < j \le l$ be integers. Let $f(s,i,j)$ denote the set of all
strings of the form $s_1 \ldots s_ics_j \ldots s_l$, where $c \in A$ is
any character that does not occur in the substring $s_i \ldots s_j$. Note
that $f(s,i,j)$ may be empty if there is no such character $c$. If
$L$ is an $l$-set and $1 \le i < j \le l$, let 
$$f(L,i,j) = \bigcup_{s \in L}\ f(s,i,j).$$

We say that the set $L'$ can be derived from the $l$-set $L$ and denote it
by $L \rightarrow L'$ if there exist integers $i,j$, $1 \le i < j \le l$ such 
that $L' = f(L,i,j)$. In the automaton, there is a transition from $L$ to $L'$
labeled $\{i,j\}$. There are ${l \choose 2}$ transitions from each $l$-set $L$,
corresponding to all possible choices of $i,j$. $L'$ is defined for each choice, 
though different choices of $i,j$ may yield the same set $L'$, including possibly the empty set.
Note that $L'$ is an $(l+i+2-j)$-set and we can derive other sets from $L'$. Let 
$\Rightarrow$ denote the transitive closure of $\rightarrow$. Thus $L \Rightarrow L'$ 
iff there exists a sequence of subsets of $A^*$ in $\mathcal{L}$, $L_1, L_2, \ldots, L_n$ 
such that $L = L_1$, $L' = L_n$ and $L_i \rightarrow L_{i+1}$, for $1 \le i < n$.

Let $S$ be the 3-set containing the string \texttt{acb}. 

We can now state the equivalence with the four color theorem.

\begin{thm}
Every planar graph is 4-colorable iff $S \not\Rightarrow \emptyset$.
\end{thm}

It is well-known that the four color theorem is true if
it is true for 4-connected plane triangulations. Whitney's theorem~\cite{W}
implies that such triangulations have a Hamiltonian cycle. Some of the
reformulations, as in ~\cite{CRZ,M}, are obtained by viewing such a triangulation
as the union of two maximal outerplanar graphs that have the edges of the
Hamiltonian cycle in common. 

Here, we view these differently. A near-triangulation is a planar 2-connected graph
in which every face except possibly the external face is a triangle. The following
lemma gives a property of 4-connected triangulations that we will use.

\begin{lem}
The vertices of any 4-connected plane triangulation $G$ can be ordered
$v_1,v_2,\ldots,v_n$ such that the subgraph $G_i$ of $G$ induced by
$\{v_1,\ldots,v_i\}$ and the subgraph $\overline{G_i}$ of $G$ induced
by $\{v_{i+1},\ldots,v_n\}$ are both near-triangulations, for
all $3 \le i \le n-3$. Also $v_1, v_2, v_n$ can be chosen to be the vertices
in the external face of a plane embedding of $G$.
\end{lem}

This property has been used elsewhere, for example in~\cite{D}, but we include the proof for completeness. 
Fix any plane embedding of the graph $G$ and let $v_1v_2$ be any edge in the boundary of the external face. Let 
$v_n$ be the third vertex in the external face of $G$. Let $v_3$ be the internal vertex such that $v_1,v_2,v_3$ 
is a face of $G$. If there is a cutvertex $v$ in $G-\{v_1,v_2,v_3\}$, since $G$ is 4-connected, 
$v_1,v_2,v_3$ must be adjacent to at least one vertex in each component of $G-\{v_1,v_2,v_3,v\}$,
contradicting the fact that $v_1,v_2,v_3$ is a face of $G$. Thus $G-\{v_1,v_2,v_3\}$ is
a near-triangulation. 

Assume that for some $i$, $3 \le i < n-3$, we have found vertices $v_1,\ldots,v_i$
such that $G_j$ and $\overline{G_j}$ are near-triangulations for $3 \le j \le i$.
Let $v_n = w_1,w_2,\ldots,w_l$ be the vertices in the boundary of the external face of
$\overline{G_i}$. If there is no chord in $\overline{G_i}$ joining two non-consecutive vertices 
$w_a, w_b$, choose $v_{i+1}$ to be any vertex $w \neq v_n$ in the external face of $\overline{G_i}$ 
that is adjacent to at least two vertices in $G_i$.  If there is a chord $w_aw_b$, let $a,b$ be such that 
$a < b$ and $b-a$ is minimum among all possible choices. Let $v_{i+1}$ be any vertex $w_j$ for $a < j < b$ 
that is adjacent to at least two vertices in $G_i$. There must exist at least one such vertex, otherwise
$w_a$ and $w_b$ have a common neighbor in $G_i$ and $G$ has a separating triangle.
This process can be continued as long as $i < n-3$. Once $v_{n-3}$ has been chosen we can choose $v_{n-2}$
and $v_{n-1}$ arbitrarily.

Note that at every step in this process, $v_{i+1}$ is adjacent to at least 2 vertices in $G_i$ and the
neighbors of $v_{i+1}$ in $G_i$ form a consecutive sequence of vertices in the boundary of the external
face of $G_i$. 

The connection between 4-coloring and strings is now clear from this. For $3 \le i \le n$, let
$v_1=w_1,w_2,\ldots,w_l = v_2$ be the vertices in the external boundary of $G_i$. Let $L_i$
be the set of all strings $g(w_1)g(w_2)\ldots g(w_l)$, where $g$ is any proper 4-coloring of $G_i$ with
colors $\{\texttt{a},\texttt{b},\texttt{c},\texttt{d}\}$. Without loss of generality, we can assume $g(v_1) = $ \texttt{a},
$g(v_2) = $ \texttt{b}, and $g(v_3) = $ \texttt{c}. So for $i=3$, the set $L_3$ contains only the
string \texttt{acb} and $L_3 = S$. If $v_{i+1}$ is adjacent to the vertices $w_j,w_{j+1},\ldots,w_k$, 
for $1 \le j < k \le l$, then the external face of $G_{i+1}$ is $w_1,\ldots,w_j,v_{i+1},w_k,\ldots,w_l$. 
The set $L_{i+1}$ of strings obtained from proper 4-colorings of $G_{i+1}$ is exactly the set 
$f(L_i,j,k)$, by definition. Thus if $S \not\Rightarrow \emptyset$, there exists a proper 4-coloring of $G$. 
On the other hand, if $S \Rightarrow \emptyset$, we can construct a near-triangulation that is not 4-colorable 
from a sequence of derivations $S = L_3 \rightarrow L_4 \rightarrow \cdots \rightarrow \emptyset$, 
using the labels of the transitions from $L_i$ to $L_{i+1}$.

A possible approach to proving this may be to identify some property of the sets $L$ such
that $S \Rightarrow L$ and show that the empty set does not satisfy it.  One such property
that follows from the  four color theorem is that any such $L$ must contain a string in which either
the character \texttt{c} or \texttt{d} does not occur. However, to prove this by induction, we
need to show that some other kinds of strings also appear in each such set. Alternatively,
characterize $l$-sets $L$ such that $L \Rightarrow \emptyset$, and show that $S$ does not
satisfy the property. A starting point may be to prove the five color theorem using this
approach with a five letter alphabet.

\end{document}